\documentclass[12pt,oneside]{article}
\usepackage{amssymb}
\author{Alexander Gorodnik}
\date{}
\title{On Oppenheim-type conjecture for systems of quadratic forms\footnote{This article is a part of author's PhD thesis at Ohio State University done under supervision of Prof.~Bergelson.}}
\newtheorem{thm}{Theorem}
\newtheorem{lem}[thm]{Lemma}

\newtheorem{cor}[thm]{Corollary}
\newtheorem{con}[thm]{Conjecture}
\newenvironment{proof}[1][Proof]{\begin{trivlist}
  \item[\hskip \labelsep {\bfseries #1}]}{\hfill$\square$\end{trivlist}}
\newenvironment{remark}[1][Remark]{\begin{trivlist}
   \item[\hskip \labelsep {\bfseries #1}]}{\end{trivlist}}
\def\qed{\end{proof}}
\def\colon{{:}\; }

\begin{document}
\maketitle
%\vspace{-1.3cm}
%\begin{center}
%{\small gorodnik@math.ohio-state.edu\\ Department of Mathematics\\ Ohio State University\\ Columbus, OH 43210\\}
%\end{center}
%\vspace{0.2cm}
\begin{abstract}
Let $Q_i$, $i=1,\ldots,t$, be real nondegenerate indefinite quad\-ratic forms in $d$ variables.
We investigate under what conditions the closure of the set $\{(Q_1(\bar x),\ldots,Q_t(\bar x))\colon\bar x\in\mathbb{Z}^d-\{\bar 0\}\}$
contains $(0,\ldots,0)$.
As a corollary, we deduce 
several results on the magnitude of the set $\Delta$ of $g\in\hbox{GL}(d,\mathbb{R})$
such that the closure of the set $\{(Q_1(g\bar x),\ldots,Q_t(g\bar x)):\;\bar x\in\mathbb{Z}^d-\{\bar 0\}\}$
contains $(0,\ldots,0)$.
Special cases are described when
depending on the mutual position of the hypersurfaces $\{Q_i=0\}$, $i=1,\ldots,t$,
the set $\Delta$ has full Haar measure or measure zero and Hausdorff dimension
$d^2-\frac{d-2}{2}$.
\end{abstract}

\section{Introduction}

It was conjectured by Oppenheim that if 
$$
Q(\bar x)=\sum_{i,j} a_{ij} x_ix_j,\;\; \bar x=(x_1,\ldots,x_d),
$$
is a real nondegenerate indefinite
quadratic form in $d$ variables, $d\ge 5$, which is not proportional to a rational form,
then
\begin{equation}\label{o1}
\forall\varepsilon>0\;\;\exists \bar x\in\mathbb{Z}^d-\{\bar 0\}\colon\; |Q(\bar x)|<\varepsilon.
\end{equation}
This conjecture was proved by Margulis \cite{mar89}.
In fact, he proved that if $Q$ is a real nondegenerate indefinite
quadratic form in $d$ variables, $d\ge 3$, which is not proportional to a rational form, then
$$
\forall\varepsilon>0\;\;\exists \bar x\in\mathbb{Z}^d-\{\bar 0\}\colon\; 0<|Q(\bar x)|<\varepsilon.
$$
See \cite{mar97} for an up-to-date survey.
It is also known that (\ref{o1}) fails for 
some indefinite irrational forms in $2$ variables.
For example, (\ref{o1}) fails for $Q=x_1^2-\alpha^2x_2^2$ with a badly approximable number $\alpha$
(e.g. $\alpha=1+\sqrt{2}$).
Recall that a number $\alpha$ is called {\bf badly approximable}
if there exists $c>0$ such that for all integers $p$ and $q>0$,
$$
\left|q\alpha -p\right|\ge \frac{c}{q}.
$$

In this paper we study similar questions for systems $(Q_j\colon j=1,\ldots,t)$ of quadratic forms.
Namely, we would like to know when the following statement is true:
\begin{equation}\label{o2}
\forall\varepsilon>0\;\;\exists \bar x\in\mathbb{Z}^d-\{\bar 0\}\colon\;\; \mathop{\max}_{j=1,\ldots, t} |Q_j(\bar x)|<\varepsilon.
\end{equation}
We expect that an analog of the Oppenheim conjecture holds for systems of quadratic forms
(see Conjecture \ref{q2_con} below).
However, the proof of this conjecture seems to be beyond the reach
of available methods. 

One can give a convenient characterization of property (\ref{o2}) for systems 
quadra\-tic forms in two variables:
(\ref{o2}) holds for nondegenerate indefinite quadratic forms
$$
Q_j(x_1,x_2)=(a_jx_1+b_jx_2)(c_jx_1+d_jx_2),\;\; j=1,\ldots, t,
$$
iff the set
$$
\bigcap_{j=1}^t\left\{\frac{a_j}{b_j}, \frac{c_j}{d_j}\right\}
$$
contains a number which is not badly approximable (we consider $\infty$ to be not badly approximable).
In the case of one quadratic form, this was observed by Dani \cite{dan00}.
Theorem \ref{thm_q2_1} below is an analog of this fact in dimension $d>2$.

Validity of (\ref{o2}) depends on the common position of the hypersurfaces $\{Q_j=0\}$, $j=1,\ldots, t$.
For example, it is easy to see that
if these hypersurfaces intersect only at the origin, then (\ref{o2}) fails.
We investigate the case when the intersection of the hypersurfaces 
$\{Q_j=0\}$, $j=1,\ldots, t$, is still relatively small: they intersect transversally in a strong sense.
To formulate our result, we recall that 
a vector $\bar\alpha=(\alpha_1,\ldots,\alpha_d)\in\mathbb{R}^d$ 
is called {\bf well approximable} of order $r$ ($r>0)$ if there exist $\bar p_n\in \mathbb{Z}^d$
and $q_n>0$ such that
\begin{equation}\label{eq_wav}
q_n^r\cdot\|q_n\bar\alpha-\bar p_n\|\rightarrow 0\;\;\hbox{as}\;\; n\rightarrow\infty.
\end{equation}

\begin{thm}\label{thm_q2_1}
Let $Q_j$, $j=1,\ldots,t$, be nondegenerate indefinite quad\-ratic forms in $d$ variables, $d\ge 3$.
Suppose that the intersection of the hypersurfaces $\{Q_j=0\}$, $j=1,\ldots, t$, is not $\{\bar 0\}$,
and at every point of intersection, which is different from $\bar 0$, the space spanned by the normal vectors
to the hypersurfaces has dimension $d-1$.
Then (\ref{o2}) holds iff the intersection of the hypersurfaces $\{Q_j=0\}$, $j=1,\ldots, t$,
contains a vector
$\bar v\in\mathbb{R}^d$
such that for some $j=1,\ldots,t$, $v_j\ne 0$, and the vector
$\left(\frac{v_i}{v_j}\colon i\ne j\right)\in\mathbb{R}^{d-1}$
is well approximable of order one.
\end{thm}

We also study the case when zero hypersurfaces of the quadratic forms have a common tangent plane at a
point of intersection.
Choose and fix a line along which the hypersurfaces $\{Q_j=0\}$, $j=1,\ldots, t$, have a common tangent plane. 
Let $\bar f_1$ be the direction vector of this line.
We choose a basis $\{\bar f_i\colon i=1,\ldots,d\}$ such that the vector $\bar f_2$ is outside
of the tangent plane, and $\bar f_i$, $i=3,\ldots,d$, are in the tangent plane. 
Define $a_t\in\hbox{\rm SL}(d,\mathbb{R})$ by
\begin{equation}\label{eq_a_t0}
a_t \bar f_1=e^{-t}\bar f_1,\quad a_t \bar f_2=e^{t}\bar f_2,\quad a_t \bar f_i=\bar f_i,\; i=3,\ldots,d.
\end{equation}
The following theorem establishes a connection between properties of the semiorbit of $a_t$, $t>0$,
in $\hbox{\rm SL}(d,\mathbb{R})/\hbox{\rm SL}(d,\mathbb{Z})$ and Diophantine condition (\ref{o2}).
First results of this type were discovered by Dani \cite{da85,da86}.

\begin{thm}\label{thm_q2_2}
Let $Q_j$, $j=1,\ldots,t$, be nondegenerate indefinite quad\-ratic forms in $d$ variables, $d\ge 3$.
Suppose that the hypersurfaces $\{Q_j=0\}$, $j=1,\ldots,t$,  have a common tangent plane at a point
of intersection, which is different from $\bar 0$. 
If the semiorbit 
$$
\{a_t\hbox{\rm SL}(d,\mathbb{Z})\colon t>0\}\subseteq \hbox{\rm SL}(d,\mathbb{R})/\hbox{\rm SL}(d,\mathbb{Z}),
$$
where $a_t$ is defined in (\ref{eq_a_t0}),
is not relatively compact, then (\ref{o2}) holds.
Moreover, if for some $\alpha_j\in \mathbb{R}$, $j=1,\ldots,t$, the quadratic form
$\sum_j\alpha_j Q_j$ is definite of rank $d-1$, then (\ref{o2}) implies that the semiorbit
$\{a_t\hbox{\rm SL}(d,\mathbb{Z})\colon t>0\}$ is not relatively compact in
$\hbox{\rm SL}(d,\mathbb{R})/\hbox{\rm SL}(d,\mathbb{Z})$.
\end{thm}

\begin{remark}
Since the choice of the vectors $\bar f_i$, $i=2,\ldots,d$, is not unique,
the transformation $a_t$ is not uniquely defined.
Nonetheless, one can check that the property that the positive semiorbit of $a_t$ 
is bounded is independent of a chosen basis.
\end{remark}

Theorems \ref{thm_q2_1} and \ref{thm_q2_2} allow us to investigate how often (\ref{o2})
holds in the space of systems of quadratic forms.
For quadratic forms $Q_j$, $j=1,\ldots, t$ in $d$ variables, we consider a parametric family
$$
(Q_j^g\colon j=1,\ldots,t),\;\; g\in \hbox{\rm GL}(d,\mathbb{R}),
$$
where $Q_j^g(\bar x)\stackrel{def}{=}Q(g\bar x)$, and denote by 
$\Delta(Q_1,\ldots,Q_t)$ the set of $g \in \hbox{\rm GL}(d,\mathbb{R})$ such that
$(Q_j^g\colon j=1,\ldots,t)$ satisfies (\ref{o2}).

\begin{cor} \label{th_m}
Let $Q_j$, $j=1,\ldots,t$, be nondegenerate indefinite quad\-ratic forms in $d$ variables, $d\ge 3$,
and $\Delta=\Delta(Q_1,\ldots,Q_t)$.
\begin{enumerate}
\item[(i)] If the intersection of the hypersurfaces $\{Q_j=0\}$, $j=1,\ldots, t$, is not $\{\bar 0\}$,
and at every point of intersection, which is different from $\bar 0$, the space spanned by the normal vectors
to the hypersurfaces has dimension $d-1$, then $\Delta$ has measure zero and Hausdorff dimension $d^2-\frac{d-2}{2}$.
\item[(ii)] If the hypersurfaces $\{Q_j=0\}$, $j=1,\ldots,t$,  have a common tangent plane at a point
of intersection, which is different from $\bar 0$, then $\Delta$ has full Haar measure in $\hbox{\rm GL}(d,\mathbb{R})$
(i.e. its complement has measure zero).
Moreover, if for some $\alpha_j\in \mathbb{R}$, $j=1,\ldots,t$, the quadratic form
$\sum_j\alpha_j Q_j$ is definite of rank $d-1$, then the complement of $\Delta$ in $\hbox{\rm GL}(d,\mathbb{R})$
has Hausdorff dimension $d^2$.
\end{enumerate}
\end{cor}

Recall that the {\bf radical} $\hbox{\rm Rad}(Q)$ of a quadratic form $Q$, which is defined on a vector space $S$, is the subspace of vectors in $S$ that are orthogonal to $S$ with respect to $Q$.

We prove the following result that complements Corollary \ref{th_m}(ii):

\begin{thm}\label{th_r}
Let $Q_j$, $j=1,\ldots,t$, be nondegenerate indefinite quad\-ratic forms in $d$ variables, $d\ge 3$,
and $\Delta=\Delta(Q_1,\ldots,Q_t)$.
Suppose that for some $\beta_j\in\mathbb{R}$, $j=2,\ldots, d$,
$$
V\stackrel{def}{=}\bigcap_{j=2,\ldots,d} \hbox{\rm Rad}(Q_j-\beta_jQ_1)
$$
has dimension at least $2$, and $V$ contains a vector $\bar v\ne \bar 0$ such that $Q_1(\bar v)=0$.
\begin{enumerate}
\item[(a)] Unless $\dim V=2$ and the quadratic form $Q_1|_V$ is nondegenerate, the complement of $\Delta$
in $\hbox{\rm GL}(d,\mathbb{R})$
is contained in a countable union of proper submanifolds of dimension at most $d^2-d+1$.
\item[(b)] If $\dim V=2$, the quadratic form $Q_1|_V$ is nondegenerate, and for some $\alpha_j\in \mathbb{R}$, $j=1,\ldots,t$, the quadratic form
$\sum_j\alpha_j Q_j$ is definite of rank $d-2$, then the complement of $\Delta$
in $\hbox{\rm GL}(d,\mathbb{R})$ has Hausdorff dimension $d^2$.
\end{enumerate}
\end{thm}

Note that the conditions of Theorem \ref{th_r} imply that the hypersurfaces $\{Q_j=0\}$, $j=1,\ldots,t$, have a common tangent plane at $\bar v$. In particular, it follows from Corollary \ref{th_m}(ii)
that the complement of $\Delta$ has measure $0$.

The simplest example for Theorem \ref{th_r}(a) is provided by the system of quadratic forms
$$
Q_j(\bar x)=x_1^2+x_2x_3+a_jx_3^2,\quad j=1,\ldots,t.
$$
It follows from the result of Dani and Margulis \cite{dm90} that the complement of
$\Delta(Q_1,Q_2)$ is a countable union of proper submanifolds of $\hbox{\rm GL}(3,\mathbb{R})$.

It is easy to refine Corollary \ref{th_m} and Theorem \ref{th_r} as follows.
Denote by $\Delta'$ the set of $g \in \hbox{\rm GL}(d,\mathbb{R})$ such that
\begin{equation}\label{eq_dp}
\forall\varepsilon>0\;\;\exists \bar x\in\mathbb{Z}^d-\{\bar 0\}\colon \;\; 0<\mathop{\max}_{j=1,\ldots, t}|Q_j(g\bar x)|<\varepsilon.
\end{equation}

\begin{cor}\label{q2_col}
Statements of Corollary \ref{th_m} and Theorem \ref{th_r} are valid for $\Delta'$.
\end{cor}

The next section is devoted to the proofs of stated results.
Some open problems are discussed in Section 3.

{\sc Acknowledgments:}
Thanks are due to my advisor V.~Bergelson for many stimulating and enlightening discussions.
I am also grateful to the referee for several pertinent suggestions, which have 
significantly improved the quality of this paper.

\section{Proofs}

Proofs of Theorems \ref{thm_q2_1} and \ref{thm_q2_2} involve only elementary linear algebra. 
Theorem \ref{thm_q2_1} follows from Lemmas \ref{q2_l3} and \ref{q2_l31} below.
To prove Theorem \ref{thm_q2_2}, we reduce it to a question about
a system of linear forms.
The proof of Corollary \ref{th_m}(i) uses Lemmas \ref{q2_l3} and \ref{q2_l31} and the result of Jarn\'\i k \cite{ja}.
Corollary \ref{th_m}(ii) and Theorem \ref{th_r} are deduced
from known properties of certain flows on the space of unimodular lattices.

If it is not stated otherwise, $\|\cdot\|$ denotes the Euclidean norm on $\mathbb{R}^d$.

The following three lemmas are used in the proof of Theorem \ref{thm_q2_1}.

\begin{lem} \label{q2_l2}
Let $Q_j$, $j=1,\ldots,t$, be nondegenerate indefinite quad\-ratic forms in $d$ variables, $d\ge 3$,
whose zero hypersurfaces intersect along a line $\ell$, and the span of normal vectors
to the hypersurfaces $\{Q_j=0\}$, $j=1,\ldots,t$, at a point of $\ell$ has dimension $d-1$.
Denote by $\pi$ the orthogonal projection on the line $\ell$.
For $m$, $M$, $\varepsilon>0$, put
\begin{eqnarray*}
R_{m,M}&=&\{\bar x\in\mathbb{R}^d\colon  \|\pi(\bar x)\|\ge M, \|\bar x-\pi(\bar x)\|\le m\|\pi(\bar x)\|\},\\
T(\varepsilon)&=&\{\bar x\in R_{m,M}\colon  \|\pi(\bar x)\|\cdot \|\bar x-\pi(\bar x)\|\le \varepsilon\},\\
S(\varepsilon)&=&\{\bar x\in R_{m,M}\colon  |Q_j(\bar x)|\le \varepsilon,j=1,\ldots,t\}.
\end{eqnarray*}
Then there exist $m,M,c_1,c_2>0$ such that for every $\varepsilon\in (0,1)$,
$$
T(c_1\varepsilon)\subseteq S(\varepsilon)\subseteq T(c_2\varepsilon).
$$
\end{lem}

\begin{proof} 
It is enough to prove the lemma for a system of quadratic forms 
$Q_j$, $j=1,\ldots,d-1$, such that normal vectors to $\{Q_j=0\}$, $j=1,\ldots,d-1$, at 
a point of $\ell$ are linearly independent.

For a vector $\bar x=(x_1,\ldots,x_d)\in\mathbb{R}^d$, put $x=(x_1,\ldots,x_{d-1})$.

Using a change of variables, we can transform $\ell$ to the line $\{x_i=0\colon i=1,\ldots,d-1\}$.
Thus, it is enough to consider the case when the line $\ell$ is as above. 
Then
\begin{equation}\label{eq_Q_i}
Q_i(\bar x)=q_i(x)+x_dL_i(x),\;\; i=1,\ldots,d-1.
\end{equation}
Note that $L_i$ is independent of $x_d$ because $Q_i(0,\ldots,0,1)=0$.
The tangent plane to $\{Q_i=0\}$ at $(0,\ldots,0,1)$ is $\{L_i=0\}$.

Let
$$
T(a,\varepsilon)=T(\varepsilon)\cap\{x_d=a\}\quad\hbox{and}\quad S(a,\varepsilon)=S(\varepsilon)\cap\{x_d=a\}.
$$
For $|a|\ge M$, we have
\begin{eqnarray}
T(a,\varepsilon)=\{(x,a)\colon  |a|\cdot \|x\|\le \varepsilon, \|x\|\le m|a|\}, \label{eq_T}
\end{eqnarray}
\begin{equation}
S(a,\varepsilon)=\left\{(x,a)\colon |q_i(x)+aL_i(x)|<\varepsilon, i=1,\ldots, d-1;\|x\|\le m |a|\right\}. \label{eq_S}
\end{equation}

First, we claim that for fixed sufficiently small $m>0$,
\begin{equation} \label{q2_eq31}
\hbox{diam} (S(a,\varepsilon))\rightarrow 0
\end{equation}
uniformly on $\varepsilon \in (0,1)$ as $|a|\rightarrow\infty$.
Let $L=(l_{ij})$ be a linear map where $l_{ij}$ are the coefficients of the linear forms $L_i$.
Note that $L$ is nondegenerate because the normal vectors to $\{Q_j=0\}$, $j=1,\ldots,d-1$,
at a point of $\ell$ are linearly independent.
Using (\ref{eq_Q_i}) and (\ref{eq_S}), we have
\begin{equation} \label{q2_eq5}
|aL_i(x)|\le\varepsilon+|q_i(x)|\le 1+\beta_i\Vert x\Vert^2\le 1+\beta_i m|a|\cdot \|x\|.
\end{equation}
for $x\in S(a,\varepsilon)$ and $i=1,\ldots, d-1$.

Suppose that there exists a sequence $({u}^{(n)},a_n)\in S(a_n,\varepsilon)$ with $|a_n|\ge 1$,
and $\Vert {u}^{(n)}\Vert\rightarrow\infty$.
We may assume that $\frac{{u}^{(n)}}{\Vert {u}^{(n)}\Vert}\rightarrow {v}_0$
for some ${v}_0$ on the unit sphere $S^{d-2}$.
From (\ref{q2_eq5}) we get 
$$
|L_i({u}^{(n)})|\le\frac{1}{|a_n|}+\beta_i m\Vert {u}^{(n)}\Vert\le 1+\beta_i m\Vert{u}^{(n)}\Vert,\;\; i=1,\ldots, d-1.
$$
It follows that 
\begin{equation}\label{eq_contr1}
|L_i({v}_0)|\le \beta_i m,\;\; i=1,\ldots,d-1.
\end{equation}
Take 
\begin{equation}\label{eq_contr2}
m=\frac{1}{2}\min\left\{\max\left\{\beta_i^{-1}|L_i({x})|\colon i=1,\ldots,d-1\right\}\colon {x}\in S^{d-2}\right\}.
\end{equation}
Note that $m>0$ because $L_i$, $i=1,\ldots, d-1$, are linearly independent.
Then the equations (\ref{eq_contr1}) and (\ref{eq_contr2}) contradict each other. 
This shows that $\hbox{diam}(S(a,\varepsilon))$ is bounded  for $\varepsilon\in (0,1)$ and $a\ge 1$.
From the first part of (\ref{q2_eq5}), we see that
$|a|\cdot \Vert L(x)\Vert$ for ${x}\in S(a,\varepsilon)$ is bounded too.
Since $L$ is nondegenerate, for some $c>0$, $\hbox{diam}(S(a,\varepsilon))<\frac{c}{|a|}$
when $\varepsilon\in (0,1)$ and $a\ge 1$. 
This proves the claim (\ref{q2_eq31}).

Consider a transformation $u\colon \mathbb{R}^{d-1}\rightarrow \mathbb{R}^{d-1}$ defined by
$$
u_i=q_i(x)+aL_i(x),\quad i=1,\ldots,d-1.
$$
By (\ref{eq_S}),
\begin{equation}\label{eq_S1}
S(a,\varepsilon)=\left\{(x,a)\colon |u_i(x)|<\varepsilon, i=1,\ldots,d-1;\|x\|\le m |a|\right\}.
\end{equation}
The Jacobian of the transformation $u$ is $J=\left(\frac{\partial q_i}{\partial x_j}\right)+aL$.
Note that $\det J({0})=\det L\ne 0$. 
Therefore by (\ref{q2_eq31}),
$u$ is a diffeomorphism on $S(a,\varepsilon)$ for sufficiently large $|a|$.
Let $g=dx_1^2+\cdots+dx_{d-1}^2$ be the standard Riemann metric,
$u^*g=du_1^2+\cdots+du_{d-1}^2$ be the pull-back of $g$
under the transformation $u$, and $L^*g$ be the pull-back of $g$ under $L$.
Clearly, for some $r_1,r_2>0$,
$$
r_1 g< L^*g< r_2g.
$$
Since $a^{-1}J\rightarrow L$ as $|a|\rightarrow\infty$
uniformly on compact sets of $\mathbb{R}^{d-1}$, $\frac{u^*g}{a^2}\rightarrow L^*g$
uniformly on compact sets. 
Therefore,
\begin{equation}\label{eq_g}
r_1g\le\frac{u^*g}{a^2}\le r_2g
\end{equation}
on $S(a,\varepsilon)$ for $|a|$ sufficiently large (say $|a|\ge M$).
Note that $\|x\|$ is the distance to the origin with respect to the metric $g$,
and $\|u(x)\|$ is the distance to the origin with respect to $u^*g$.
Therefore, it follows from (\ref{eq_g}) that for $x\in S(a,\varepsilon)$ and $|a|\ge M$,
$$
r_1^{1/2}|a|\cdot\|x\|\le\|u(x)\|\le r_2^{1/2}|a|\cdot\|x\|.
$$
By (\ref{eq_T}) and (\ref{eq_S1}), for $|a|\ge M$,
$$
T\left(a,r_2^{-1/2}\varepsilon\right)\subseteq S(a,\varepsilon)\subseteq T\left(a,r_1^{-1/2}\sqrt{2}\varepsilon\right).
$$
The lemma is proved.
\qed

\begin{lem} \label{q2_l3}
For a line $\ell$ in $\mathbb{R}^d$, denote by $\pi_\ell$ the orthogonal projection on the line $\ell$.
Assume that zero hypersurfaces of a system of nondegenerate indefinite quadratic forms $Q_j$, $j=1,\ldots,t$,
in $d$ variables intersect at a point, which is different from $\bar 0$, and at every 
nonzero point of intersection, the span of the normal vectors to these hypersurfaces has dimension $d-1$.
Then  (\ref{o2}) holds iff there exist a line $\ell$
in the intersection of the hypersurfaces of $\{Q_j=0\}$, $j=1,\ldots,t$, and a sequence $\bar x^{(n)}\in\mathbb{Z}^d-\{\bar 0\}$ such that
\begin{equation}\label{eq_cond}
\Vert\pi_\ell (\bar{x}^{(n)})\Vert\rightarrow\infty\quad\hbox{and}\quad\Vert\pi_\ell (\bar{x}^{(n)})\Vert\cdot \|\bar{x}^{(n)}-\pi_\ell(\bar{x}^{(n)})\|\rightarrow 0.
\end{equation}
\end{lem}

\begin{proof}  Let $m,M,c_1,c_2>0$ be as in Lemma \ref{q2_l2}.

Suppose that $Q_j$, $j=1,\ldots,t$, satisfy (\ref{o2}).
Then there exists a sequence $\bar{x}^{(n)}\in\mathbb{Z}^d-\{\bar 0\}$ such that 
\begin{equation}\label{eq_en}
\varepsilon_n\stackrel{def}{=}\max\{ |Q_j(\bar{x}^{(n)})|\colon j=1,\ldots, t \}\rightarrow 0.
\end{equation}
Passing if needed to a subsequence, we may assume that 
$$
\bar{y}^{(n)}\stackrel{def}{=}\frac{\bar{x}^{(n)}}{\Vert\bar{x}^{(n)}\Vert}\rightarrow \bar{y}_0
$$
for some $\bar{y}_0$
on the unit sphere. 
Then $Q_j(\bar{y}_0)=0$, $j=1,\ldots,t$. 
Let $\ell$ be the line through $\bar{y}_0$ and the origin.

If $\Vert \bar{x}^{(n)}\Vert\nrightarrow\infty$, we can take $\{\bar{x}^{(n)}\}$ to be bounded. 
Then for each $j$, $Q_j(\bar{x}^{(n)})$ admits only finitely many values for $n\ge 1$.
Therefore, by (\ref{eq_en}),
$Q_j(\bar{x}^{(n)})=0$, $j=1,\ldots,t$, for sufficiently large $n$.
We can replace $\bar{x}^{(n)}$ by $n\bar{x}^{(n)}$.
Thus, we may assume that $\Vert \bar{x}^{(n)}\Vert\rightarrow\infty$.

Since $\bar{y}^{(n)}\rightarrow\bar{y}_0$, the sequence $\bar{x}^{(n)}$ is inside of the cone
$$
\left\Vert\bar{x}-\pi_\ell(\bar{x})\right\Vert < m \Vert\pi_\ell(\bar{x})\Vert
$$
for sufficiently large $n$.
It follows that if $\|\pi_\ell(\bar{x}^{(n)})\|\nrightarrow\infty$, then $\Vert \bar{x}^{(n)}\Vert\nrightarrow\infty$.
Therefore, $\pi_\ell(\bar{x}^{(n)})\rightarrow\infty$. 
In particular, $\bar{x}^{(n)}\in R_{m,M}$ for sufficiently
large $n$. 
By Lemma \ref{q2_l2},
$$
\Vert\pi_\ell (\bar{x}^{(n)})\Vert\cdot \|\bar{x}^{(n)}-\pi_\ell(\bar{x}^{(n)})\|\le c_2\varepsilon_n\rightarrow 0. 
$$

Conversely, suppose that (\ref{eq_cond}) is satisfied for some line $\ell$ in the intersection of
the hypersurfaces $\{Q_j=0\}$, $j=1,\ldots, t$.
Then for sufficiently large $n$, $\bar{x}^{(n)}\in R_{m,M}$.
Applying Lemma \ref{q2_l2}, we get
$$
|Q_j(\bar{x}^{(n)})|\le c_1^{-1} \Vert\pi_\ell (\bar{x}^{(n)})\Vert\cdot \|\bar{x}^{(n)}-\pi_\ell(\bar{x}^{(n)})\|\rightarrow 0,\;\; j=1,\ldots,t.
$$
This proves the lemma.
\qed

\begin{lem} \label{q2_l31}
Let $\bar v=(v_1,\ldots,v_{d})\in\mathbb{R}^d$ with 
$v_j\ne 0$ for some $j=1,\ldots, d$, and $\pi$ is the orthogonal projection on the direction $\bar v$.
Then the vector $\left(\frac{v_i}{v_j}\colon i\ne j\right)\in\mathbb{R}^{d-1}$
is well approximable of order $1$ iff for some sequence $\bar x^{(n)}\in\mathbb{Z}^d-\{\bar 0\}$,
\begin{equation}\label{eq_cond2}
\Vert\pi (\bar{x}^{(n)})\Vert\rightarrow\infty\quad\hbox{and}\quad\Vert\pi (\bar{x}^{(n)})\Vert\cdot \|\bar{x}^{(n)}-\pi(\bar{x}^{(n)})\|\rightarrow 0.
\end{equation}
\end{lem}

\begin{proof} 
Without loss of generality, $j=1$.
Put $\alpha_i=\frac{v_{i+1}}{v_1}$, $i=1,\ldots, d-1$.

Let $L_i(\bar x)=\alpha_i x_1-x_{i+1}$, $i=1,\ldots, d-1$.
Note that the planes
defined by $L_i=0$, $i=1,\ldots,d-1$, intersect along the line in direction $\bar v$.

Suppose that $(\alpha_1,\ldots,\alpha_{d-1})$ is well approximable of order $1$. 
There exists a sequence
$\bar{x}^{(n)}\in\mathbb{Z}^d$ such that
\begin{equation}\label{eq_1}
|x_1^{(n)}|\cdot \left|L_i\left(\bar x^{(n)}\right)\right|\rightarrow 0,\;\; i=1,\ldots, d-1,
\end{equation}
with $x_1^{(n)}\ne 0$ for all $n$. 
It follows that
\begin{equation}\label{eq_li}
L_i(\bar x^{(n)})\rightarrow 0,\;\; i=1,\ldots, d-1.
\end{equation}

Suppose that $|x_1^{(n)}|$ is bounded.
A linear map $A\colon \mathbb{R}^d\rightarrow\mathbb{R}^d$, which is defined by
\begin{equation}\label{eq_AA}
(A\bar x)_i=L_i(\bar x),\; i=1,\ldots,d-1,\quad\hbox{and}\quad (A\bar x)_d=x_1,
\end{equation}
is nondegenerate. 
Hence, by (\ref{eq_li}), 
$\bar x^{(n)}$ is bounded, and $L_i(\bar x^{n})=0$, $i=1,\ldots, d-1$,
for sufficiently large $n$. 
In this case, (\ref{eq_cond2}) holds for the sequence $n\bar x^{(n)}$.

If $|\bar x_1^{(n)}|$ is not bounded, we may assume that $|\bar x_1^{(n)}|\rightarrow\infty$.
By (\ref{eq_li}), for sufficiently large $n$,
\begin{equation}\label{eq_B}
\|A\bar x^{(n)}\|=\left(\sum_{i=1}^{d-1}L_i(\bar x^{(n)})^2+(x_1^{(n)})^2\right)^{1/2}\le \sqrt{d} |x_1^{(n)}|.
\end{equation}
Since $A$ is nondegenerate, for some $c_1>0$ and every $\bar x\in\mathbb{R}^d$, 
\begin{equation}\label{eq_pi}
\|\pi(\bar x)\|\le c_1 \|A\bar x\|.
\end{equation}
Let $P$ be the plane through the origin
orthogonal to $\bar v$. 
A map $C\colon P\rightarrow\mathbb{R}^{d-1}$
defined by
$$
(C\bar x)_i=L_i(\bar x),\;\; i=1,\ldots,d-1,
$$
is invertible. 
Therefore, for some
$c_2>0$ and every $\bar x\in P$,
\begin{equation}\label{eq_B2}
\|\bar x-\pi(\bar x)\|\le c_2\max\{|L_i(\bar x)|\colon i=1,\ldots,d-1\}.
\end{equation}
In fact, the last inequality holds for every $\bar x\in\mathbb{R}^d$ because it is independent
of translations by vectors parallel to $\bar v$. 
Now (\ref{eq_cond2}) follows from 
(\ref{eq_1}), (\ref{eq_B}), (\ref{eq_pi}), (\ref{eq_B2}).

Conversely, suppose that (\ref{eq_cond2}) holds for some sequence $\bar x^{(n)}\in \mathbb{Z}^d-\{\bar 0\}$.
Note that for $\bar x\in\mathbb{R}^d$, $|L_i(\bar x)|$ is up to a constant the distance from $\bar x$
to the plane $\{L_i=0\}$, and $\|\bar x-\pi(\bar x)\|$ is the distance from $\bar x$ to the line 
through $\bar v$.
Therefore, for some $r_i>0$ and every $\bar x\in\mathbb{R}^d$,
\begin{equation}\label{eq_di}
|L_i(\bar x)|\le r_i\|\bar x-\pi(\bar x)\|,\;\; i=1,\ldots, d-1.
\end{equation}
In particular, it follows that 
\begin{equation}\label{eq_li2}
L_i(\bar x^{(n)})\rightarrow 0,\;\; i=1,\ldots,d-1.
\end{equation}
Clearly, $x_1^{(n)}\ne 0$ for sufficiently large $n$.

Consider a linear map $B$ defined by
$$
(B\bar x)_i=L_i(\bar x),\; i=1,\ldots,d-1,\quad\hbox{and}\quad (B\bar x)_d=\pi(\bar x).
$$
Since $B$ is nondegenerate, for some $r_3>0$ and every $\bar x\in \mathbb{R}^d$,
\begin{equation}\label{eq_A1}
|x_1|\le r_3 \|B\bar x\|.
\end{equation}
By (\ref{eq_cond2}) and (\ref{eq_li2}), for sufficiently large $n$,
\begin{equation}\label{eq_A}
\|B\bar x^{(n)}\|=\left(\sum_{i=1}^{d-1}L_i(\bar x^{(n)})^2+\|\pi(\bar x^{(n)})\|^2 \right)^{1/2}
\le \sqrt{d}\|\pi(\bar x^{(n)})\|.
\end{equation}
Finally, (\ref{eq_1}) follows from (\ref{eq_cond2}), (\ref{eq_di}), (\ref{eq_A1}), and (\ref{eq_A}).
Thus, the vector $(\alpha_1,\ldots,\alpha_{d-1})$ is well approximable of order $1$.
\qed

Combining Lemmas \ref{q2_l3} and \ref{q2_l31}, we deduce Theorem \ref{thm_q2_1}.

The proof of Corollary \ref{th_m}(i) uses Lemma \ref{q2_l3}, Lemma \ref{q2_l31},
and the result of Jarn\'\i k on
Hausdorff dimension of the set of well approximable vectors \cite{ja}.
It is known that the intersection of the set of well approximable
vectors of order $r$ in $\mathbb{R}^d$ with every nonempty open set has Hausdorff dimension $\frac{d+1}{r+1}$.

\begin{proof} [Proof of Corollary \ref{th_m}(i)]
Let $\mathcal{Q}$ be the projective variety in the complex projective space $\mathbb{P}^{d-1}(\mathbb{C})$
defined by $Q_j=0$, $j=1,\ldots, t$. 
Since at every point of $\mathcal{Q}(\mathbb{R})$
the rank of Jacobian $\frac{\partial(Q_1,\ldots,Q_t)}{\partial(x_1,\ldots,x_d)}$
is equal to $d-1$, it follows that every irreducible component of $\mathcal{Q}$ that 
has nonempty intersection with $\mathcal{Q}(\mathbb{R})$ has dimension zero.
Therefore, $\mathcal{Q}(\mathbb{R})$ consists of finitely many points.
This means that the zero hypersurfaces $\{Q_j=0\}$, $j=1,\ldots, t$, 
in $\mathbb{R}^d$ intersect along only finitely many lines
that pass through the origin. 
Let $\bar v_s$, $s=1,\ldots, N$, be the direction
vectors of these lines. 
Then zero hypersurfaces of $Q_j^g$, $j=1,\ldots, t$,
intersect along vectors $g^{-1}\bar v_s$, $s=1,\ldots, N$.
For some $g_s\in \hbox{\rm GL}(d,\mathbb{R})$, $\bar v_s=g_s\bar e_1$, where $e_1=(1,0,\ldots,0)$. 
We have
$$
g^{-1}\bar v_s=\left( (g^{-1}g_s)_{11}, \ldots,(g^{-1}g_s)_{d1}\right),\;\; s=1,\ldots,N.
$$
By Lemma \ref{q2_l3}, 
$$
\Delta=\bigcup_{s=1}^N \Delta_s,
$$
where $\Delta_s$ is the set of $g\in\hbox{\rm GL}(d,\mathbb{R})$ such that (\ref{eq_cond})
holds for the line $\ell$ through the origin in direction $g^{-1}\bar v_s$. 
Clearly, it is
enough to prove the theorem for each of the sets $\Delta_s$ separately.

Fix $s=1,\ldots,N$.
Let
$$
G_{s,i}=\{g\in\hbox{\rm GL}(d,\mathbb{R})\colon (g^{-1}g_s)_{i1}\ne 0 \},
$$
and $\Delta_{s,i}=\Delta_s\cap G_{s,i}$ for $i=1,\ldots, d$.

By Lemma \ref{q2_l31}, for $g\in G_{s,i}$,
$g$ belongs to $\Delta_{s,i}$ iff the vector 
$$
\left( \frac{(g^{-1}g_s)_{j1}}{(g^{-1}g_s)_{i1}}\colon j\ne i\right)
$$
is well approximable of order $1$.
Consider maps $F_{s,i}\colon G_{s,i}\rightarrow \hbox{\rm M}(d,\mathbb{R})$,  $i=1,\ldots,d$,
defined by
\begin{eqnarray*}
F_{s,i}(g)_{j1}&=&\frac{(g^{-1}g_s)_{j1}}{(g^{-1}g_s)_{i1}}\;\;\;\hbox{for}\;\; j\ne i,\\
F_{s,i}(g)_{j_1 j_2}&=&(g^{-1}g_s)_{j_1 j_2}\;\; \hbox{for}\;\;(j_1,j_2)\notin \{(j,1),\; j\ne i\}.
\end{eqnarray*}
Note that $F_{s,i}$ is a diffeomorphism on $G_{s,i}$, and $F_{s,i}(G_{s,i})$ is a complement of a finite
union of hypersurfaces in $\hbox{\rm M}(d,\mathbb{R})$.
Let $W_i$ be the set of $g\in \hbox{\rm M}(d,\mathbb{R})$ such that the vector 
$(g_{j1}\colon j\ne i)$ is well-approximable of order $1$.
We have that for $g\in G_{s,i}$, $g$ belongs to $\Delta_{s,i}$ iff $F_{s,i}(g)$ belongs to $W_i$, i.e.
$\Delta_{s,i}=F_{s,i}^{-1}\left(W_i\cap F_{s,i}(G_{s,i})\right)$. 
Therefore,
\begin{equation}\label{eq_d}
\Delta_s=\bigcup_{i=1}^{d}\Delta_{s,i}=\bigcup_{i=1}^{d} F_{s,i}^{-1}\left(W_i\cap F_{s,i}(G_{s,i})\right).
\end{equation}

The intersection of the set of well approximable vectors of order $1$ in $\mathbb{R}^{d-1}$
with every nonempty open subset has Hausdorff dimension $\frac{d}{2}$.
It follows that the intersection of the set $W_i$, $i=1,\ldots, d$, with every nonempty open subset
of $\hbox{\rm M}(d,\mathbb{R})$
has Hausdorff dimension $d(d-1)+1+\frac{d}{2}=d^2-\frac{d-2}{2}$
(see \cite[Section~3.5.5]{bd}).
Since $F_{s,i}(G_{s,i})$ is open in $\hbox{\rm M}(d,\mathbb{R})$,
the set $W_i\cap F_{s,i}(G_{s,i})$ has the same property.
Therefore, by (\ref{eq_d}), the intersection of $\Delta_s$ with every nonempty open subset of
$\hbox{\rm M}(d,\mathbb{R})$ has Hausdorff dimension $d^2-\frac{d-2}{2}$.
\qed

The following two lemmas are crucial for the proofs of Theorem \ref{thm_q2_2} and 
Corollary \ref{th_m}(ii). 
The proofs of the lemmas are based on ideas of Dani \cite{da85,da86}.

\begin{lem} \label{l_dani}
Let $L_i$, $i=1,\ldots, d$, be linearly independent linear forms in $d$ variables, $d\ge 2$.
Define $a_t\in\hbox{\rm SL}(d,\mathbb{R})$ by
\begin{eqnarray}\label{eq_a_t}
L_1(a_t\bar x)&=&e^{-t}L_1(\bar x),\\
L_2(a_t\bar x)&=&e^{t}L_2(\bar x),\\
L_i(a_t\bar x)&=&L_i(\bar x),\;\; i=3,\ldots,d.\nonumber
\end{eqnarray}
Then the set $\{a_t\hbox{\rm SL}(d,\mathbb{Z})\colon t>0\}$ is not bounded in $\hbox{\rm SL}(d,\mathbb{R})/\hbox{\rm SL}(d,\mathbb{Z})$
iff there exists a sequence $\bar x^{(n)}\in\mathbb{Z}^d-\{\bar 0\}$ such that
\begin{equation}\label{eq_dani}
L_1(\bar x^{(n)})L_2(\bar x^{(n)})\rightarrow 0\;\;\hbox{and}\;\; L_i(\bar x^{(n)})\rightarrow 0,\;\; i=2,\ldots, d.
\end{equation}
\end{lem}

\begin{proof} 
Let $\|\cdot\|$ be the norm on $\mathbb{R}^d$ defined by
$$
\|\bar x\|=\max\{|L_i(\bar x)|\colon i=1,\ldots,d\},\quad \bar x\in\mathbb{R}^d.
$$
By Mahler Compactness Criterion, the orbit $\{a_t\hbox{\rm SL}(d,\mathbb{Z})\colon t>0\}$, is unbounded iff
for some $t_n>0$ and $\bar x^{(n)}\in\mathbb{Z}^d-\{\bar 0\}$,
$a_{t_n}\bar x^{(n)}\rightarrow \bar 0$, i.e.
\begin{eqnarray}\label{eq_dani1}
\|a_{t_n}\bar x^{(n)}\|=\max\{e^{-t_n}|L_1(\bar x^{(n)})|, e^{t_n}|L_2(\bar x^{(n)})|,&\\
|L_i(g\bar x^{(n)})|\colon i=3,\ldots,d\}
&\rightarrow 0.
\nonumber
\end{eqnarray}
It is clear that (\ref{eq_dani1}) implies (\ref{eq_dani}).
Conversely, suppose that (\ref{eq_dani}) holds.

First, we consider the case when $L_2(\bar x^{(n)})=0$ for infinitely many $n$. 
Taking a subsequence, we may assume that $L_2(\bar x^{(n)})=0$ for all $n$. 
Then (\ref{eq_dani1}) holds for every sequence $t_n\rightarrow \infty$ such that
$$
e^{-t_n}|L_1(\bar x^{(n)})|\rightarrow 0.
$$

Now we may assume that $L_2(\bar x^{(n)})\ne 0$.
If $L_1(\bar x^{(n)})$ is bounded, then since the linear forms $L_i$, $i=1,\ldots,d$, are linearly independent,
$\bar x^{(n)}$ is bounded. 
This implies that for sufficiently large $n$, $L_2(\bar x^{(n)})=0$.
Thus, we may assume that $|L_1(\bar x^{(n)})|\rightarrow\infty$.
Then (\ref{eq_dani1}) holds for the sequence $t_n$ such that 
$e^{t_n}=|L_1(\bar x^{(n)})|^{1/2}\cdot|L_2(\bar x^{(n)})|^{-1/2}$. 
Note that $t_n>0$ for sufficiently large $n$. 
This proves the lemma.
\qed

\begin{lem} \label{l_dani00}
Let $L_i$, $i=1,\ldots, d$, be linearly independent linear forms in $d$ variables, $d\ge 2$.
Then the set $\Phi$ of $g\in\hbox{\rm GL}(d,\mathbb{R})$ such that 
\begin{equation}\label{eq_LLL2}
L_1(g\bar x^{(n)})L_2(g\bar x^{(n)})\rightarrow 0\;\;\hbox{and}\;\; L_i(g\bar x^{(n)})\rightarrow 0,\;\; i=2,\ldots, d.
\end{equation}
for some sequence $\bar x^{(n)}\in\mathbb{Z}^d-\{\bar 0\}$ has full Haar measure, and its
complement has Hausdorff dimension $d^2$.
\end{lem}

\begin{proof} 
Let $L_i^0(\bar x)=x_i$, $i=1,\ldots, d$. 
For some $g_0\in\hbox{\rm GL}(d,\mathbb{R})$,
$L_i(g\bar x)=L_i^0(g_0g\bar x)$, $i=1,\ldots, d$. 
Since multiplication by $g_0$ is a diffeomorphism of $\hbox{\rm GL}(d,\mathbb{R})$,
we may assume without loss of generality that $g_0=1$ and $L_i=L_i^0$, $i=1,\ldots,d$.

For every matrix $u\in \hbox{\rm GL}(d,\mathbb{R})$,
$$
u=\hbox{diag}(\lambda,1,\ldots,1)\cdot v
$$
with $\lambda\in\mathbb{R}^\times$ and $v\in\hbox{\rm SL}(d,\mathbb{R})$.
Moreover, the map
$$
\hbox{\rm GL}(d,\mathbb{R})\to\mathbb{R}^\times\times\hbox{\rm SL}(d,\mathbb{R})\colon  u\mapsto (\lambda,v)
$$
is a diffeomorphism.
Since
$$
L_1(u\bar x)=\lambda L_1(v\bar x),\quad L_i(u x)=L_i(v\bar x),\;\; i=2,\ldots, d,
$$
(\ref{eq_LLL2}) holds for $g=u$ iff it holds for $g=v$.
Therefore, if we show that the set $\Phi\cap\hbox{\rm SL}(d,\mathbb{R})$
has full measure in $\hbox{\rm SL}(d,\mathbb{R})$,
this will imply that the set $\Phi$ has full measure in $\hbox{\rm GL}(d,\mathbb{R})$ too.
Also by the argument as in \cite[Section~3.5.5]{bd}, 
if the set $\hbox{\rm SL}(d,\mathbb{R})-(\Phi\cap\hbox{\rm SL}(d,\mathbb{R}))$
has Hausdorff dimension $d^2-1$, then the set $\hbox{\rm GL}(d,\mathbb{R})-\Phi$
has Hausdorff dimension $d^2$. 
Hence, the proof reduces to a question about $\hbox{\rm SL}(d,\mathbb{R})$.

By Lemma \ref{l_dani}, (\ref{eq_LLL2}) holds for some sequence $\bar x^{(n)}\in\mathbb{Z}^d-\{\bar 0\}$
iff the set 
$$
\{g^{-1}a_tg\hbox{\rm SL}(d,\mathbb{Z})\colon t>0\}
$$
is unbounded in $\hbox{\rm SL}(d,\mathbb{R})/\hbox{\rm SL}(d,\mathbb{Z})$.
This is equivalent to the set $\{a_tg\hbox{\rm SL}(d,\mathbb{Z})\colon t>0\}$
being unbounded.
It is known that the set of $g\in\hbox{\rm SL}(d,\mathbb{R})$ such that the orbit
$\{a_tg\hbox{\rm SL}(d,\mathbb{Z})\colon t>0\}$ is bounded has Haar measure zero
(Moore Ergodicity Theorem \cite[Theorem 2.2.6]{zim84}),
and its intersection with every nonempty open subset has Hausdorff
dimension $d^2-1$ (Kleinbock, Margulis \cite{km}).
Therefore, the set of $g\in\hbox{\rm SL}(d,\mathbb{R})$
such that (\ref{eq_LLL2}) holds for some sequence $\bar x^{(n)}\in\mathbb{Z}^d-\{\bar 0\}$
has full measure, and its complement in $\hbox{\rm SL}(d,\mathbb{R})$ has
Hausdorff dimension $d^2-1$. 
This implies the lemma.
\qed

\begin{proof} [Proof of Theorem \ref{thm_q2_2}]
With respect to the basis $\{\bar f_i\colon  i=1,\ldots,d\}$,
\begin{equation}\label{eq_Q_j}
Q_j(\bar y)= q_j(y_2,\ldots,y_{d})+y_1l_j(y_2,\ldots,y_{d}),\;\; j=1,\ldots,t,
\end{equation}
where $q_j$ is a quadratic form, and $l_j$ is a linear form. 
Note that $l_j$ is
independent of $y_1$ because $Q_j(\bar f_1)=0$. 
The tangent plane to $\{Q_j=0\}$ at $\bar f_1$ is $\{l_j=0\}$. 
Thus, $l_j=c_j l$, $j=1,\ldots,t$, for some $c_j\in\mathbb{R}-\{0\}$
and a linear form $l$. 
Since $\bar f_1$ and $\bar f_i$, $i=3,\ldots,d$, are in the tangent plane,
$l(\bar y)=y_2$. 
Let $L_i(\bar x)$, $i=1,\ldots, d$, denote the coordinates of vector $\bar x$
with respect to the basis $\{\bar f_i\colon i=1,\ldots,d\}$. 
Then
$$
Q_j(\bar x)= q_j(L_2(\bar x),\ldots,L_{d}(\bar x))+c_jL_1(\bar x)L_2(\bar x),\;\; j=1,\ldots,t,
$$
and (\ref{eq_a_t}) holds.
It follows that (\ref{o2}) holds
for $(Q_j\colon j=1,\ldots,t)$ provided that for some sequence $\bar x^{(n)}\in\mathbb{Z}^d-\{\bar 0\}$,
\begin{equation}\label{eq_LLL}
L_1(\bar x^{(n)})L_2(\bar x^{(n)})\rightarrow 0\quad\hbox{and}\quad L_i(\bar x^{(n)})\rightarrow 0,\; i=2,\ldots, d.
\end{equation}
Thus, the first statement of the theorem follows from Lemma \ref{l_dani}.

Suppose that $Q=\sum_j \alpha_j Q_j$ is definite of rank $d-1$. 
Since $Q$ is definite,
it follows from (\ref{eq_Q_j}) that $Q(\bar y)=Q(y_2,\ldots,y_d)$. 
Therefore, (\ref{o2}) implies (\ref{eq_LLL})
for some sequence $\bar x^{(n)}\in\mathbb{Z}^d-\{\bar 0\}$,
and the second statement of the theorem follows from Lemma \ref{l_dani}.
\qed

\begin{proof} [Proof of Corollary \ref{th_m}(ii)]
It was shown in the proof of Theorem \ref{thm_q2_2}
that (\ref{o2}) holds for $(Q_j^g\colon j=1,\ldots,t)$ provided that for some sequence
$\bar x^{(n)}\in\mathbb{Z}^d-\{\bar 0\}$,
\begin{equation}\label{eq_dani000}
L_1(g\bar x^{(n)})L_2(g\bar x^{(n)})\rightarrow 0\quad\hbox{and}\quad L_i(g\bar x^{(n)})\rightarrow 0,\; i=2,\ldots, d,
\end{equation}
and when a linear combination of $Q_j$, $j=1,\ldots,d$, is definite of rank $d-1$, 
$(Q_j^g\colon j=1,\ldots,t)$ satisfies (\ref{o2}) iff (\ref{eq_dani000}) holds.
Thus, Corollary \ref{th_m}(ii) follows from Lemma \ref{l_dani}.
\qed

To prove Theorem \ref{th_r}(a), we use the following lemma:
\begin{lem}\label{l_unip}
Let $Q_j$, $j=1,\ldots,t$, be a system of quadratic forms that satisfies the conditions
of Theorem \ref{th_r}(a). 
Then there exists a nontrivial one-parameter unipotent subgroup of $\hbox{\rm SL}(d,\mathbb{R})$
which leaves all $Q_j$, $j=1,\ldots,t$, invariant.
\end{lem}

\begin{proof} 
Let $Q_1|_V$ be nondegenerate. 
Then $\dim V\ge 3$. 
Write 
$$
\mathbb{R}^d=V\oplus W
$$
where $W$ is the orthogonal complement of
$V$ with respect to $Q_1$. 
Then the group 
$$
H\stackrel{def}{=}\hbox{SO}(Q_1|_V)\oplus id_W
$$
leaves $Q_1$
invariant. 
Since $V$ is in the radical of $Q_j-\beta_j Q_1$, the forms $Q_j$, $j=2,\ldots,t$
are invariant under $H$ too. 
The group $H$ is a noncompact semisimple
Lie group. 
Hence, it contains a nontrivial one-parameter unipotent subgroup. 
This proves the lemma when $Q_1|_V$ is nondegenerate. 

Now suppose that $Q_1|_V$ is degenerate, and $Q_1|_V\ne 0$.
Choose a basis $\{\bar f_i\colon i=1,\ldots,d\}$ of $\mathbb{R}^d$
such that $\bar f_1$ is in the radical of $Q_1|_V$,
$\bar f_2$ is in $V$, $Q_1(\bar f_2)\ne 0$, and $\bar f_i$, $i=3,\ldots, t$, are orthogonal to
$\bar f_2$ with respect to $Q_1$. 
Since $\bar f_2$ is in the radical of $Q_j-\beta_j Q_1$,
$j=2,\ldots,t$, $\bar f_2$ is orthogonal to $\bar f_i$, $i=3,\ldots,t$, with respect
to $Q_j$, $j=2,\ldots,t$. 
In the basis $\{\bar f_i\colon i=1,\ldots,d\}$,
$$
Q_j(\bar y)=\beta_j (y_2^2+y_1L(y_3,\ldots,y_d))+q_j(y_3,\ldots,y_d),\;\; j=1,\ldots,t
$$
for some linear form $L$ and quadratic forms $q_j$, $j=1,\ldots,t$. 
(Here we put $\beta_1=1$.)
Note that $L$ is the same all $j$ because $\bar f_1$ is in the radical of
$Q_j-\beta_j Q_1$, $j=2,\ldots,t$.
Define a linear transformation $u_t\in\hbox{\rm SL}(d,\mathbb{R})$ by
\begin{eqnarray*}
u_t&:&y_1\mapsto y_1-2ty_2-t^2L(y_3,\ldots, y_d),\\
u_t&:&y_2\mapsto y_2+tL(y_3,\ldots, y_d),\\
u_t&:&y_i\mapsto y_i,\;\; i\ge 3.
\end{eqnarray*}
Then $\{u_t\colon t\in\mathbb{R}\}$ is a nontrivial one-parameter unipotent group that stabilizes $Q_j$, $j=1,\ldots,t$.

It remains to consider the case when $Q_1|_V=0$.
Choose a basis $\{\bar f_i\colon i=1,\ldots,d\}$, such that $\bar f_1$ and $\bar f_2$ is
in $V$. 
With respect to this basis,
$$
Q_j(\bar y)=\beta_j (y_1L_1(y_3,\ldots,y_d)+y_2L_2(y_3,\ldots,y_d))+q_j(y_3,\ldots,y_d)
$$
for some linear forms $L_1$, $L_2$, and quadratic forms $q_j$, $j=1,\ldots,t$.
(Here we put $\beta_1=1$.) 
Note that $L_1$ and $L_2$ 
are independent of $j$ because $\bar f_1$ and $\bar f_2$ are
in the radical of $Q_j-\beta_j Q_1$, $j=2,\ldots,t$.
The linear forms $L_1$ and $L_2$ are not zero because the
quadratic form $Q_1$ is nondegenerate.
Define a linear transformation $v_t\in\hbox{\rm SL}(d,\mathbb{R})$ by
\begin{eqnarray*}
v_t&:&y_1\mapsto y_1+tL_2(y_3,\ldots, y_d),\\
v_t&:&y_2\mapsto y_2-tL_1(y_3,\ldots, y_d),\\
v_t&:&y_i\mapsto y_i,\;\; i\ge 3.
\end{eqnarray*}
Then $\{v_t\colon t\in\mathbb{R}\}$ is a nontrivial one-parameter unipotent group that stabilizes $Q_j$, $j=1,\ldots,t$.
We have proved the lemma.
\qed

The following Lemma is used in the proof of Theorem \ref{th_r}(b).
Its proof is essentially the same as the proof of Lemma \ref{l_dani} and is omitted.

\begin{lem} \label{l_dani2}
Let $L_i$, $i=1,\ldots, d$, be linearly independent linear forms in $d$ variables, $d\ge 2$,
and $a_t$ is defined as in (\ref{eq_a_t}).
Then the set $\{a_t\hbox{\rm SL}(d,\mathbb{Z})\colon t\in\mathbb{R}\}$ is not bounded in $\hbox{\rm SL}(d,\mathbb{R})/\hbox{\rm SL}(d,\mathbb{Z})$
iff there exists a sequence $\bar x^{(n)}\in\mathbb{Z}^d-\{\bar 0\}$ such that
\begin{equation}\label{eq_dani_new}
L_1(\bar x^{(n)})L_2(\bar x^{(n)})\rightarrow 0\;\;\hbox{and}\;\; L_i(\bar x^{(n)})\rightarrow 0,\;\; i=3,\ldots, d.
\end{equation}
\end{lem}

\begin{proof} [Proof of Theorem \ref{th_r}]
Let $G=\hbox{\rm SL}(d,\mathbb{R})$ and $\Gamma=\hbox{\rm SL}(d,\mathbb{Z})$.

To prove (a), we show first that the complement of $\Delta\cap G$ in $G$
is contained in a countable union of submanifolds of dimension at most $d^2-d$.
Denote by $U\subset G$ a one-parameter unipotent subgroup
that stabilizes all $Q_j$, $j=1,\ldots,t$. 
Such a subgroup exists by Lemma \ref{l_unip}.
By Ratner's topological rigidity \cite{ra}, the set $\Omega$ of $g\in G$
such that $Ug\Gamma$ is not dense in $G/\Gamma$
is a countable union of sets of the form $F\Gamma$ where
$F$ is a connected Lie subgroup of $G$.
It is known that the maximal dimension of a proper connected subgroup
of $G$ is $d^2-d$ \cite[Sec.~3.3]{ov}.
By Mahler compactness criterion, for every $g\in G-\Omega$, 
there exist $u_n\in U$ and $\bar x^{(n)}\in\mathbb{Z}^d-\{\bar 0\}$ such that $u_ng\bar x^{(n)}\to \bar 0$ as
$n\to\infty$. 
Then
$$
Q_j(g\bar x^{(n)})=Q_j(u_ng\bar x^{(n)})\to 0,\;\; j=1,\ldots,t.
$$
Therefore, $G-\Delta\cap G$ is contained in $\Omega$.

Let $\hbox{\rm GL}^\pm(d,\mathbb{R})$ be the group of matrices with positive/negative discriminant. 
It is clear that
\begin{equation}\label{eq_glp}
\Delta\cap\hbox{\rm GL}^+(d,\mathbb{R})=\bigcup_{\lambda>0} \lambda \cdot(\Delta\cap G)
\end{equation}
Therefore, the complement of $\Delta\cap\hbox{\rm GL}^+(d,\mathbb{R})$ in $\hbox{\rm GL}^+(d,\mathbb{R})$
is contained in a countable union of submanifolds of dimension at most $d^2-d+1$.
Similarly,
\begin{equation}\label{eq_glm}
\Delta\cap\hbox{\rm GL}^-(d,\mathbb{R})=\bigcup_{\lambda>0} \lambda g_0\cdot(\tilde\Delta\cap G)
\end{equation}
for a fixed $g_0\in\hbox{\rm GL}^-(d,\mathbb{R})$ and $\tilde\Delta\stackrel{def}{=}\Delta(Q_1^{g_0},\ldots,Q_t^{g_0})$.
It follows that the complement of $\Delta\cap\hbox{\rm GL}^-(d,\mathbb{R})$ in $\hbox{\rm GL}^-(d,\mathbb{R})$
is contained in a countable union of submanifolds of dimension at most $d^2-d+1$.
This proves (a).

Now we prove (b).
Take a basis $\{\bar f_i\colon i=1,\ldots,d\}$ such
that $V=\left<\bar f_1,\bar f_2\right>$, and $\bar f_i$, $i=3,\ldots,d$ is orthogonal to
$V$ with respect to $Q_1$. 
It follows from the definition of $V$ that
$\bar f_i$, $i=3,\ldots,d$ are orthogonal to $V$ with respect to $Q_j$,
$j=2,\ldots,t$, too. 
In addition, we can change $\bar f_1$ and $\bar f_2$ such that
with respect to the basis $\{\bar f_i\colon i=1,\ldots,d\}$,
\begin{equation}\label{eq_qb}
Q_j(\bar y)=\beta_jy_1y_2+q_j(y_3,\ldots,y_d),\;\; j=1,\ldots,t
\end{equation}
for some quadratic forms $q_j$, $j=1,\ldots,t$. 
(Here we put $\beta_1=1$.)
Let $Q\stackrel{def}{=}\sum_j \alpha_j Q_j$ be a definite form of 
rank $d-2$. 
It follows from (\ref{eq_qb}) that $Q(\bar y)=Q(y_3,\ldots,y_d)$.
Denote by $L_i(\bar x)$, $i=1,\ldots, d$, the coordinates of a vector $\bar x$
with respect to the basis $\{\bar f_i\colon i=1,\ldots,d\}$. 
Then (\ref{o2}) holds for
$(Q_j^g\colon j=1,\ldots,t)$ with $g\in G$
iff for some sequence $\bar x^{(n)}\in\mathbb{Z}^d-\{\bar 0\}$,
$$
L_1(g\bar x^{(n)})L_2(g\bar x^{(n)})\rightarrow 0\;\;\hbox{and}\;\; L_i(g\bar x^{(n)})\rightarrow 0,\;\; i=3,\ldots, d.
$$
By Lemma \ref{l_dani2}, this is equivalent to the set
$\{g^{-1}a_tg\Gamma\colon t\in\mathbb{R}\}$
being unbounded in $G/\Gamma$. 
Therefore, by the result of Kleinbock and Margulis \cite{km},
the set $G-(\Delta\cap G)$ has Hausdorff dimension $d^2-1$.

To show that the complement of $\Delta$ in $\hbox{\rm GL}(d,\mathbb{R})$ has Hausdorff dimension $d^2$,
one may use (\ref{eq_glp}), (\ref{eq_glm}), and \cite[Section~3.5.5]{bd}. 
This proves (b).
\qed

\begin{proof} [Proof of Corollary \ref{q2_col}]
Under the conditions in (i), the intersection of hypersurfaces $Q_j=0$, $j=1,\ldots,t$,
consists of finitely many lines that pass through the origin. 
Denote by
$\bar v_s$, $s=1,\ldots,N$, the direction vectors of these lines.
Note that for $g\in\hbox{\rm GL}(d,\mathbb{R})$, hypersurfaces $Q_j^g=0$, $j=1,\ldots,t$,
intersect along lines in directions $g^{-1}\bar v_s$, $s=1,\ldots,N$.
For $s=1,\ldots,N$ and $\bar x\in\mathbb{R}^d-\{\bar 0\}$, denote by 
$R_s(\bar x)$ the set of $g\in\hbox{\rm GL}(d,\mathbb{R})$ such that
$g^{-1}\bar v_s$ is collinear with $\bar x$. 
The set $R_s(\bar x)$ is a submanifold
of dimension $d^2-(d-1)$. 
Since
$$
\Delta-\Delta'\subseteq\bigcup_{1\le s\le N;\;\bar x\in \mathbb{Z}^d-\{\bar 0\}} R_s(\bar x),
$$
the set $\Delta-\Delta'$ has Hausdorff dimension at most $d^2-(d-1)<d^2-\frac{d-2}{2}$. 
This shows
that the set $\Delta'$ has Hausdorff dimension $d^2-\frac{d-2}{2}$.

To prove (ii), we define
$$
K_j(\bar x)=\{g\in\hbox{\rm M}(d,\mathbb{R})\colon  Q_j(g\bar x)=0\},\;\; j=1,\ldots, t,\;\bar x\in \mathbb{R}^d.
$$
This is an algebraic subset of $\hbox{\rm M}(d,\mathbb{R})$, and for $\bar x\ne 0$,
$K_j(\bar x)\ne \hbox{\rm M}(d,\mathbb{R})$, which implies that it has measure $0$.
Note that 
\begin{equation}\label{eq_last}
\Delta-\Delta'\subseteq\bigcup_{1\le j\le t;\;\bar x\in \mathbb{Z}^d-\{\bar 0\}} K_j(\bar x),
\end{equation}
where the last set has measure $0$ and Hausdorff dimension $<d^2$.
Therefore, when $\Delta$ has complement in $\hbox{\rm GL}(d,\mathbb{R})$
of measure zero, so does $\Delta'$. 
The statement about Hausdorff dimension follows from (\ref{eq_last}) too.
This proves (ii).

We use notations from the proof of Theorem \ref{th_r}.
If we show that $G-\Delta'\cap G$ is contained in $\Omega$, then 
it is possible to finish the proof as in Theorem \ref{th_r}(a).
There exists $\bar x^0\in\mathbb{R}^d$ such that $Q_j(\bar x^0)\ne 0$ for every $j=1,\ldots,t$.
Therefore, for every $\varepsilon>0$, there exists $\bar x_\varepsilon\in\mathbb{R}^d$
such that
$$
0<|Q_j(\bar x_\varepsilon)|<\varepsilon,\;\;j=1,\ldots,t.
$$
Write $\bar x_\varepsilon=g_\varepsilon \bar e_1$ for some $g_\varepsilon\in G$ and $e_1=(1,0,\ldots,0)$.
For $g\notin\Omega$, $\overline{Ug\Gamma}=G$. 
Thus, there exist $u\in U$ and $\gamma\in\Gamma$
such that
$$
0<|Q_j(ug\gamma \bar e_1)|=|Q_j(g\gamma \bar e_1)|<\varepsilon,\;\;j=1,\ldots,t.
$$
Hence, $g\in \Delta'$. 
This show that $G-\Delta'\cap G\subseteq\Omega$ and proves Theorem \ref{th_r}(a).

In the case (b), the intersection of the hypersurfaces $\{Q_j=0\}$, $j=1,\ldots,t$,
consists of two lines. 
Thus, by the same argument as in the proof of Corollary \ref{q2_col}(i),
the set $\Delta-\Delta'$ has Hausdorff dimension at most $d^2-(d-1)$. 
This implies part (b) of Theorem \ref{th_r}.
\qed

\section{Open Problems}

Theorems \ref{thm_q2_1} and \ref{thm_q2_2} 
provide information about property (\ref{o2}) only when the dimension of the space
spanned by normal vectors to the hypersurfaces $\{Q_j=0\}$, $j=1,\ldots,t$, is either
$1$ or $d-1$. 
Clearly, this leaves a gap for $d>3$.
It is of interest to investigate property (\ref{o2})
assuming other conditions on the intersection of the hypersurfaces $\{Q_j=0\}$, $j=1,\ldots,t$.
At least, one should try to give a complete answer to the question about
the magnitude of the set $\Delta(Q_1,\ldots,Q_t)$, which was partially studied in Corollary \ref{th_m}.

Corollary \ref{th_m} illustrates that in some cases, property (\ref{o2}) has
complicated Diophantine nature. 
This is analogous to the situation with Oppenheim conjecture in dimension $2$.
In higher dimensions, the following conjecture seems plausible:

\begin{con}\label{q2_con}
Let $d\in\mathbb{N}$ be sufficiently large. 
Let $(Q_1,Q_2)$ be a pair of real nondegenerate quadratic forms
in $d$ variables such that every linear combination $\alpha Q_1+\beta Q_2$ with $\alpha^2+\beta^2\ne 0$ is indefinite,
has rank $\ge 3$, and does not have all rational coefficients. 
Then (\ref{o2}) holds.
\end{con}

By a theorem of P.~Finsler, for $d\ge 3$ every linear combinations $\alpha Q_1+\beta Q_2$ 
with $\alpha^2+\beta^2\ne 0$ is
indefinite iff the intersection of the hypersurfaces $\{Q_1=0\}$ and $\{Q_2=0\}$
is not equal to $\{\bar 0\}$.
(This theorem was proved independently by many authors (see \cite{uh}).)
Therefore, under the conditions of Conjecture \ref{q2_con}, the set 
$$
\{\bar x\colon  |Q_i(\bar x)|<\varepsilon, i=1,2\}
$$
is not compact.

A question analogous to this conjecture was studied by Dani and Margulis \cite{dm90}.
They showed that if $Q$ is a nondegenerate indefinite quadratic form, and $L$ is a linear
form in $3$ variables such that the plane $\{L=0\}$ is tangent to the hypersurface $\{Q=0\}$, and every linear combination
$\alpha Q+\beta L^2$ with $\alpha^2+\beta^2\ne 0$
does not have all rational coefficients, then $\{(Q(\bar x),L(\bar x))\colon \bar x\in\mathbb{R}^d\}$
is dense in $\mathbb{R}^2$. 
Similar result holds for a pair consisting of a
linear form and a quadratic form in $d$ variables, $d\ge 4$ \cite{g1}.

Some partial results towards Conjecture \ref{q2_con} were obtained in \cite{co} by R.~J.~Cook.
He studied pairs of diagonal
quadratic forms in $d$ variables, $d\ge 9$, with algebraic coefficients.

Let us illustrate the conjecture by several examples:
\begin{enumerate}
\item Let $\mathbb{R}^d=V_1\oplus V_2$ be a direct sum of vector spaces $V_1$ and $V_2$
of dimension at least $3$. 
Let $Q_1$ and $Q_2$ be indefinite quadratic forms such that $Q_i|_{V_i}$ is nondegenerate,
$Q_i|_{V_j}=0$ for $i\ne j$, and $V_1\perp V_2$ with respect to $Q_1$ and $Q_2$.
According to Conjecture \ref{q2_con}, (\ref{o2}) should
hold for $(Q_1,Q_2)$ provided that every nonzero linear combination of $Q_1$ and $Q_2$ does not
have all rational coefficients. 
This seems to be one of the easiest special cases of the
conjecture. 
It should be possible to attack this case using the original approach
of Margulis \cite{mar89} and Ratner's topological rigidity \cite{ra}.

\item Another promising case of Conjecture \ref{q2_con} is when a pair $(Q_1,Q_2)$ of quadra\-tic forms
satisfies the conditions of Theorem \ref{th_r}(a). 
By Lemma \ref{l_unip},
both $Q_1$ and $Q_2$ are invariant under a nontrivial unipotent subgroup.
Hence, one can use Ratner's topological rigidity.

\item In general, the group of linear transformations that leaves both forms invariant
may be finite. 
If this is the case, the method of Margulis, which is based on dynamics on
homogeneous spaces of Lie groups, does not work. 
This problem can appear 
in every dimension even when the hypersurfaces $\{Q_1=0\}$ and $\{Q_2=0\}$ have a common tangent plane. 
For example, when
$$
Q_1(\bar x)=\sum_{i=1}^d x_{d-i+1}x_i,\quad
Q_2(\bar x)=\sum_{i=1}^{d-1} x_{d-i}x_i+\alpha Q_1(\bar x),
$$
up to a linear change of variables, the group that stabilizes both $Q_1$ and $Q_2$ is finite.
\end{enumerate}

%\newpage

\noindent {\sc Department of Mathematics, University of Michigan, Ann Arbor, MI 48109}\\
\noindent {\it E-mail}: \texttt{gorodnik@umich.edu}

\end{document}